\documentclass{article}

\usepackage{amsmath}
\usepackage{amssymb}
\usepackage{amsthm}
\usepackage{a4wide}
\usepackage[english]{babel}

\numberwithin{equation}{section}

\newtheorem{Th}{Theorem} 
\newtheorem{Lem}{Lemma}
\theoremstyle{remark}

\newtheorem{Rem}{\bf Remark}

\newenvironment{Proof}[1][Proof]{\proof[\rm\textbf{#1}]}{\endproof}

\newcommand{\eps}{\varepsilon}

\newcommand{\Z}{\mathbb{Z}}
\newcommand{\Q}{\mathbb{Q}}
\newcommand{\R}{\mathbb{R}}

\newcommand{\calD}{{\mathcal D}}

\newcommand{\Span}{\operatorname{Span}}
\newcommand{\Card}{\operatorname{Card}}

\renewcommand{\d}{\mathrm{d}}

\title{Many odd zeta values are irrational}
\author{St\'ephane Fischler, Johannes Sprang and Wadim Zudilin}

\date{\today}

\begin{document}

\maketitle 

\begin{abstract}
Building upon ideas of the second and third authors, we prove that at least $2^{(1-\varepsilon)\frac{\log s}{\log \log s}}$ values of the Riemann zeta function at odd integers between 3 and $s$ are irrational, where $\varepsilon$ is any positive real number and $s$ is large enough in terms of $\varepsilon$.
This lower bound is asymptotically larger than any power of $\log s$; it improves on the bound $ \frac{1-\varepsilon}{1+\log 2} \log s $ that follows from the Ball--Rivoal theorem.
 
The proof is based on construction of several linear forms in odd zeta values with related coefficients.
\end{abstract}

\section*{Introduction}

When $s\geq 2$ is an even integer, the value $\zeta(s)$ of the Riemann zeta function is a non-zero rational multiple of $\pi^s$ and, therefore, a transcendental number.  On the other hand, no such relation is expected to hold for  $\zeta(s)$ when $s\geq 3$ is odd; a folklore conjecture states that the numbers $\pi$, $\zeta(3)$, $\zeta(5)$, $\zeta(7),\ldots$ are algebraically independent over the rationals. This conjecture is predicted by Grothendieck's period conjecture for mixed Tate motives. But both conjectures are far out of reach and we do not even know the transcendence of a single odd zeta value.

\bigskip

 It was only in 1978 when Ap\'ery astonished the mathematics community by his proof \cite{Apery} of the irrationality of $\zeta(3)$ (see  \cite{SFBou} for a survey). 
 The next breakthrough was taken in 2000 by Ball and Rivoal \cite{BR, RivoalCRAS} who proved the following:

\begin{Th}[Ball--Rivoal] \label{thBR}
Let $\eps > 0$. Then for any $s\geq 3$ odd and sufficiently large with respect to $\eps$, we have
$$\dim_\Q \Span_\Q ( 1, \, \zeta(3), \, \zeta(5), \, \zeta(7), \, \ldots , \zeta(s)) \geq \frac{1-\eps}{1+\log 2}\,\log s.$$
\end{Th}

Their corresponding result for small $s$ has been refined several times \cite{Zudilincentqc, SFZu}, 
but  the question whether $\zeta(5)$ is irrational remains open.
The proof of Theorem~\ref{thBR} involves the well-poised hypergeometric series
\begin{equation}\label{eqBR}
n!^{s-2r}\, \sum_{t=1}^\infty \frac{ \prod_{j=0}^{ (2r+1)n} (t-rn+j) }{ \prod_{j=0}^{ n} (t+j)^{s+1}},
\end{equation}
which happens to be a $\Q$-linear combination of 1 and odd zeta values when $s$ is odd and $n$ is even, and Nesterenko's linear independence criterion \cite{Nesterenkocritere}. The bound $\frac{1-\eps}{1+\log 2}\log s$ follows from comparison of how small the linear combination is with respect to the size of its coefficients, after multiplying by a common denominator to make them integers. To improve on this bound using the same strategy, one has to find linear combinations that are considerably smaller, with not too large coefficients,\,---\,it comes out to be a rather difficult task. This may be viewed as
an informal explanation of why the lower bound in Theorem \ref{thBR} has never been improved for large values of $s$, whereas the theorem itself has been generalized to several other families of numbers.

\bigskip

Using (with $s=20$) the series 
$$n!^{s-6}\,  \sum_{k=1}^\infty \left.\Big( \frac{\d}{\d t}\Big)^2 \bigg( \Big(t+\frac{n}{2}\Big) \frac{\prod_{j=0}^{ 3n} (t-n+j)^3 }{ \prod_{j=0}^{ n} (t+j)^{s+3}}\bigg)\right|_{t=k} ,$$
which is a $\Q$-linear combination of 1 and odd zeta values starting from $\zeta(5)$, Rivoal has proved \cite{vingtetun} that among the numbers $\zeta(5)$, $\zeta(7)$, \ldots, $\zeta(21)$, at least one is irrational. This result has been improved by the third author \cite{Zudilinonze}: among the four numbers $\zeta(5)$, $\zeta(7)$, $\zeta(9)$, $\zeta(11)$, at least one is irrational; and he also showed \cite{Zudilincentqc} that, for any odd $\ell \geq 1$, there is an irrational number among $\zeta(\ell+2)$, $\zeta(\ell+4)$, \ldots, $\zeta(8\ell-1)$.
Proofs of these results do not require use of linear independence criteria: if a sequence of $\Z$-linear combinations of real numbers from a given (fixed) collection tends to 0, and is non-zero infinitely often, then at least one of these numbers is irrational. A drawback of this approach is that it only allows one to prove that {\em one} number in a family is irrational.

\bigskip

The situation has drastically changed when the third author introduced \cite{Zudilintrick} a new method (see also \cite{KrattZ}). He casts (with $s=25$) the rational function in the form
$$ R(t) =  2^{6n}n!^{s-5}\, \frac{\prod_{j=0}^{6n} (t-n+\frac{j}{2})}{ \prod_{j=0}^{ n} (t+j)^{s+1}}$$
and proves that both series 
$$\sum_{t=1}^\infty R(t) \quad\mbox{and}\quad \sum_{t=1}^\infty R\Big(t+\frac12\Big)$$
are $\Q$-linear combinations of 1, $\zeta(3)$, $\zeta(5)$, \ldots, $\zeta(s)$ with {\em related\/} coefficients.
This allows him to eliminate one odd zeta value, and to prove that {\em at least two} zeta values among $\zeta(3)$, $\zeta(5)$, \ldots, $\zeta(25)$ are irrational. In view of Ap\'ery's Theorem, the result means that one number among $\zeta(5)$, \ldots, $\zeta(25)$ is irrational\,---\,nothing really novel, but the method of proof is new and more elementary
than the ones in \cite{vingtetun} and \cite{Zudilinonze}
as it avoids use of the saddle point method. More importantly, the method allows to prove the irrationality of at least two zeta values in a family without having to produce very small linear forms.
The same strategy has been adopted by Rivoal and the third author \cite{RZnote} to prove that among $\zeta(5)$, $\zeta(7)$, \ldots, $\zeta(69)$, at least two numbers are irrational.
 
\bigskip

The method in \cite{Zudilintrick} has been generalized by the second author \cite{Sprang}, who introduces another integer parameter $D>1$ and considers the rational function
\begin{equation}\label{eqSprang}
R(t) = D^{6(D-1)n}n!^{s-3D-1}\,\frac{ \prod_{j=0}^{3Dn} (t-n+\frac{j}{D})}{ \prod_{j=0}^{ n} (t+j)^{s+1}}.
\end{equation}
He proves that for any divisor $d$ of $D$ the series
$$\sum_{j=1}^d \sum_{t=1}^\infty R\Big(t+\frac{j}{d}\Big)$$
is a $\Q$-linear combination of 1, $\zeta(3)$, $\zeta(5)$, \ldots, $\zeta(s)$.
The crucial point of this construction is that each $\zeta(i)$ appears in this $\Q$-linear combination with a coefficient that depends on $d$ in a very simple way. This makes it possible to eliminate from the entire collection of these linear combinations as many odd zeta values as the number of divisors of $D$. Finally, taking $D$ equal to a power of 2 and $s$ sufficiently large with respect to $D$, the second author proves that at least $\frac{\log D}{\log 2}$ numbers are irrational among $\zeta(3)$, $\zeta(5)$, \ldots, $\zeta(s)$. This strategy represents a new proof that $\zeta(i)$ is irrational for infinitely many odd integers $i$.

\bigskip

Building upon the approach in \cite{Zudilintrick} and \cite{Sprang} we prove the following result.
 
\begin{Th} \label{th1}
Let $\eps > 0$, and $s\geq 3$ be an odd integer sufficiently large with respect to $\eps$. Then among the numbers
$$\zeta(3), \, \zeta(5), \, \zeta(7), \, \ldots , \zeta(s),$$
at least
$$2^{(1-\eps)\frac{\log s}{\log \log s}}$$
are irrational.
\end{Th}

In this result, the lower bound is asymptotically greater than $\exp(\sqrt{\log s})$, and than any power of $\log s$; ``to put it roughly, [it is] much more like a power of $s$ than a power of $\log s$'' \cite[Chapter XVIII, \S 1]{HW}.

In comparison, Theorem~\ref{thBR} gives only $\frac{1-\eps}{1+ \log 2} \log s$ irrational odd zeta values, but they are linearly independent
over the rationals, whereas Theorem~\ref{th1} ends up only with their irrationality.

\bigskip

Our proof of Theorem~\ref{th1} follows the above-mentioned strategy of the second and third authors. The main new ingredient, compared to the proof in \cite{Sprang}, is taking $D$ large (about $s^{1-2\eps}$) and equal to the product of the first prime numbers (the so-called primorial)\,---\,such a number has asymptotically the largest possible number of divisors with respect to its size (see \cite[Chapter XVIII, \S 1]{HW}). To perform the required elimination of a prescribed set of odd zeta values, we need to establish
that a certain auxiliary matrix is invertible. Whereas the second author's choice of $D$ in \cite{Sprang} allows him to deal with elementary properties of a Vandermonde matrix, we use at this step a generalization of the corresponding result.
We give three different proofs of the latter, based on arguments from combinatorics of partitions, from linear algebra accompanied with a lemma of Fekete, and from analysis using Rolle's theorem.

\bigskip

The structure of this paper is as follows. In \S \ref{sec2} we construct linear forms in values of the Hurwitz zeta function. Denominators of the coefficients are studied in \S \ref{sec3}; and the asymptotics of the linear forms are dealt with in \S \ref{sec4}. Section \ref{sec5} is devoted to the proof that an auxiliary matrix is invertible. Finally, we establish Theorem~\ref{th1} in \S \ref{sec6}.

\section{Construction of linear forms} \label{sec2}

From now on we let $s$, $D$  be positive integers such that $s \geq 3D$; we assume that $s$ is odd. Let $n$ be a positive integer, such that $Dn $ is even. Consider the following rational function:
$$R_n(t) = D^{3Dn} \, \, n!^{s+1-3D} \, \, \frac{ \prod_{j=0}^{3Dn} (t-n+\frac{j}{D})}{ \prod_{j=0}^{ n} (t+j)^{s+1}}$$
which, of course, depends also on $s$ and $D$. Notice that the difference of the function from the corresponding one in \cite{Sprang} is in the factor $D^{3Dn}$ instead of $D^{6(D-1)n}$ (see Eq.~\eqref{eqSprang}).

Similar rational functions have already been considered, see \cite{Catalan} for the case $D=2$ and \cite{Nash, Nishimoto, SFcaract} for general $D$. However the ``central'' factors $t-n+\frac{j}{D}$ with $Dn < j < 2Dn$ are missing, and (as the second author noticed \cite{Sprang}) they play a central role in the arithmetic estimates (see Lemma \ref{lemarith}  below).

\bigskip

\begin{Rem} \label{Rem-r}
Though one can implement an additional parameter $r$ in the definition of the rational function $R_n(t)$,
in a way similar to the one for the Ball--Rivoal series \eqref{eqBR}, we have verified that this does not bring any
improvement to the result of Theorem~\ref{th1}.
\end{Rem}

\bigskip

The rational function $R_n(t)$ has a partial fraction expansion
\begin{equation} \label{eqelsples}
R_n(t) = \sum_{i=1}^s \sum_{k=0}^n \frac{a_{i,k}}{(t+k)^i}.
 \end{equation}
For any $j\in \{1,\ldots,D\}$, take
$$r_{n,j} = \sum_{m=1}^\infty R_n\Big(m+\frac{j}{D}\Big).$$

\bigskip

We recall that the Lerch and Hurwitz zeta functions are defined by
$$
\Phi(z,i,\alpha) = \sum_{n=0}^{\infty}\frac{z^n}{(n+\alpha)^i}
\quad\mbox{and}\quad
\zeta(i,\alpha) = \Phi(1,i,\alpha) = \sum_{n=0}^\infty \frac1{(n+\alpha)^i},
$$
where  $\alpha>0$ and also $i\geq 2$ for the latter.

\bigskip

The following is precisely \cite[Lemma~1.5]{Sprang}; the change of the normalizing factor $D^{3Dn}$ does not affect the statement.

\begin{Lem} \label{lemconstru}
For each $j\in \{1,\ldots,D\}$, we have
$$r_{n,j} 
 = \rho_{0,j}+\sum_{\substack{3\leq i \leq s\\i \;\mbox{\scriptsize odd}}} \rho_i \, \zeta\Big(i, \frac{j}{D}\Big),$$
where
$$\rho_i = \sum_{k=0}^n a_{i,k} \quad\mbox{for $3\leq i \leq s$, \ $i$ odd},$$
does not depend on $j$, and 
\begin{equation} \label{eqdefrhoz}
 \rho_{0,j} = - \sum_{k=0}^n \sum_{\ell=0}^k \sum_{i=1}^s \frac{ a_{i,k}}{(\ell+\frac{j}{D})^i} .
 \end{equation}
\end{Lem}

\begin{Proof}
We follow the strategy of proofs in \cite[Lemma~3]{Zudilintrick} and \cite[Lemma~1.5]{Sprang}. Let $z$ be a real number such that $0<z<1$. We have
\begin{align*}
\sum_{m=1}^\infty R_n\Big(m+\frac{j}{D}\Big) z^m 
&= \sum_{m=1}^\infty \sum_{i=1}^s \sum_{k=0}^n \frac{a_{i,k} z^m }{(m+k+\frac{j}{D})^i} \\
&= \sum_{i=1}^s \sum_{k=0}^n a_{i,k} z^{-k} \sum_{m=1}^\infty\frac{ z^{m+k} }{(m+k+\frac{j}{D})^i} \\
&= \sum_{i=1}^s \sum_{k=0}^n a_{i,k} z^{-k} \bigg( \Phi\Big(z, i, \frac{j}{D}\Big) - \sum_{\ell = 0}^k \frac{z^\ell}{(\ell+\frac{j}{D})^i} \bigg).
 \end{align*}
Now we let $z$ tend to 1 in the equality we have obtained; the left-hand side tends to $r_{n,j} $. On the right-hand side, the term involving the Lerch function with $i=1$ has coefficient $\sum_{k=0}^n a_{1,k} z^{-k} $. Since $ \Phi(z,1,\frac{j}{D}) $ has only a logarithmic divergence as $z\to 1$ and
$$
\sum_{k=0}^n a_{1,k} =\lim_{t\to \infty } t R_n(t) =0,
$$
this term tends to 0 as $z\to1$. All other terms have finite limits as $z\to 1$, so that
$$r_{n,j} = \rho_{0,j}+\sum_{i=2}^s \rho_i \, \zeta\Big(i, \frac{j}{D}\Big),$$
where $ \rho_{0,j} $ is given by Eq.~\eqref{eqdefrhoz}, and $ \rho_i = \sum_{k=0}^n a_{i,k} $ for any $i\in \{2,\ldots,s\}$. 
 
 \bigskip
 
To complete the proof, we apply the symmetry phenomenon of \cite{BR, RivoalCRAS}. Since $s$ is odd and $Dn $ is even we have $R_n(-n-t) = - R_n(t)$. Now the partial fraction expansion \eqref{eqelsples} is unique, so that $a_{i,n-k} = (-1)^{i+1} a_{i,k}$ for any $i$ and $k$. This implies that $\rho_i = 0$ when $i$ is even, and Lemma~\ref{lemconstru} follows.
\end{Proof}

\section{Arithmetic estimates}\label{sec3}

As usual we let $d_n = \operatorname{lcm}(1,2,\ldots,n)$.

\begin{Lem} \label{lemarith} 
We have
\begin{equation}
\label{arith-I}
d_n^{s+1-i} \rho_i \in\Z \quad\mbox{for}\; i = 3, 5, \ldots, s,
\end{equation}
and
\begin{equation}
\label{arith-II}
d_{n+1}^{s+1} \rho_{0,j} \in\Z \quad\mbox{for any}\; j\in \{1,\ldots,D\}.
\end{equation}
\end{Lem}

For part \eqref{arith-I} we use the strategy of the proof of \cite[Lemma~4.5]{SFcaract}; note that \cite[Lemma~1.3]{Sprang} does not apply in our present situation   because of the different normalization of the rational function $R_n(t)$ compared to the one in \eqref{eqSprang}.
To establish \eqref{arith-II} we follow the proof of \cite[Lemma~1.4]{Sprang};
we use $d_{n+1}$ here instead of $d_n$ to include the case corresponding to $j=D$.

\begin{Proof}[Proof of Lemma~\ref{lemarith}]
For any $\alpha\in\frac1{D}\Z$ we introduce
$$F_\alpha(t) = D^n \, \, \frac{\prod_{j=1}^n (t+\alpha+\frac{j}{D})}{\prod_{j=0}^n( t+j)} = \sum_{k=0}^n \frac{A_{\alpha,k}}{t+k},$$
where $A_{\alpha,k}$ is an integer in view of the explicit formulas
$$(-1)^k A_{\alpha,k} = \binom{n}{k} \frac{\prod_{j=1}^n (D( \alpha-k)+j )}{n!} = \begin{cases}
\binom{n}{k} \binom{D( \alpha-k) +n}{n} &\mbox{if}\; \alpha-k\geq 0, \\
0 &\mbox{if}\; \frac{-n}{D} \leq \alpha-k < 0, \\
(-1)^n \binom{n}{k} \binom{D(k- \alpha) -1}{n} &\mbox{if}\; \alpha-k < \frac{-n}{D} .
\end{cases}$$
We also consider
$$G(t) = \frac{n!}{\prod_{j=0}^n (t+j)} = \sum_{k=0}^n \frac{ (-1)^k \binom{n}{k} }{t+k},$$
so that 
\begin{equation} \label{eqpro}
R_n(t) = (t-n) \, G(t)^{s+1-3D} \, \prod_{\ell = 0 }^{3D-1} F_{-n+\frac{\ell n}{D}}(t).
\end{equation}
From this expression we compute the partial fraction expansion of $R_n(t)$ using the rules
$$ 
\frac{t-n}{t+k} = 1 - \frac{k+n}{t+k}
\quad\mbox{and}\quad
 \frac{1}{(t+k)(t+k')} = \frac{1}{(k'-k)(t+k)}+\frac{1}{(k- k')(t+k')} \quad\mbox{for}\; k\neq k'.$$
 A denominator appears each time the second rule is applied, and the denominator is always a divisor of $d_n$ (see \cite{Colmez} or \cite[Lemma~1]{Zudilintrick}). This happens $s+1-i$ times in each term that contributes to $a_{i,k}$ because there are $s+1$ factors in the product \eqref{eqpro} (apart from $t-n$).
Therefore,
 $$d_{n}^{s+1-i} a_{i,k} \in \Z \quad\mbox{for any $i$ and $k$},$$
implying \eqref{arith-I}.

\bigskip

We now proceed with the second part of Lemma~\ref{lemarith}, that is, with demonstrating the inclusions \eqref{arith-II}.
Recall from Lemma~\ref{lemconstru} that
\begin{equation} \label{eqrhoz}
d_{n+1}^{s+1} \rho_{0,j} = - \sum_{k=0}^n \sum_{\ell=0}^k \bigg( \sum_{i=1}^s \frac{d_{n+1}^{s+1} a_{i,k}}{(\ell+\frac{j}{D})^i}\bigg).
\end{equation}
If $j=D$ then 
$$d_{n+1}^{s+1-i} a_{i,k} \quad\text{and}\quad \frac{d_{n+1}^i}{(\ell+\frac{j}{D})^i}$$
are integers for any $k$, $\ell$ and $i$, so that $d_{n+1}^{s+1} \rho_{0,j} \in\Z$. 
From now on, we assume that $1\leq j \leq D-1$ and we prove that for any $k$ and any $\ell$ the internal sum over $i$ in Eq.~\eqref{eqrhoz} is an integer. With this aim in mind, fix integers $k_0$ and $\ell_0$, with $0\leq \ell_0 \leq k_0 \leq n$, and assume that the corresponding sum is not an integer. Since $1\leq j \leq D-1$ we have $R_n(\ell_0-k_0+\frac{j}{D}) = 0$, so that 
\begin{equation} \label{eqentier}
\sum_{i=1}^s \frac{d_{n+1}^{s+1} a_{i,k_0}}{(\ell_0+\frac{j}{D})^i}
= - \sum_{\substack{k=0\\k\neq k_0}}^n \sum_{i=1}^s \frac{d_{n+1}^{s+1} a_{i,k}}{(\ell_0-k_0+k +\frac{j}{D})^i} .
\end{equation}
This rational number is not an integer: it has negative $p$-adic valuation for at least one prime number $p$. Therefore, on either side of \eqref{eqentier} there is at least one term with negative $p$-adic valuation: there exist $i_0,i_1\in\{1,\ldots,s\}$ and $k_1\in \{0,\ldots,n\}$, $k_1\neq k_0$, such that 
$$v_p\bigg( \frac{d_{n+1}^{s+1} a_{i_0,k_0}}{( \ell_0+\frac{j}{D})^{i_0}} \bigg) < 0 \quad \mbox{and} \quad
v_p\bigg( \frac{d_{n+1}^{s+1} a_{i_1,k_1}}{( \ell_0-k_0+k_1+\frac{j}{D})^{i_1}} \bigg) < 0.
$$
Since $d_{n+1}^{s+1-i} a_{i,k}\in\Z$ for any $i$ and $k$, this leads to
$$v_p\bigg( \frac{d_{n+1}^{i_0} }{( \ell_0+\frac{j}{D})^{i_0}} \bigg) < 0
\quad \mbox{and} \quad
v_p\bigg( \frac{d_{n+1}^{i_1} }{( \ell_0-k_0+k_1+\frac{j}{D})^{i_1}} \bigg) < 0,
$$
implying
$$\min\bigg( v_p \Big( \ell_0+\frac{j}{D} \Big) , \, v_p \Big( \ell_0-k_0+k_1+\frac{j}{D} \Big)\bigg)> v_p(d_{n+1}).$$
As $k_0 - k_1 = ( \ell_0+\frac{j}{D} ) - ( \ell_0-k_0+k_1+\frac{j}{D} )$, we deduce that $v_p(k_0 - k_1 ) > v_p(d_{n+1})$, which is impossible in view of the inequality $0 < |k_0 - k_1 | \leq n$. The contradiction completes the proof of Lemma~\ref{lemarith}.
\end{Proof}

\bigskip

\begin{Rem} \label{remarith}
It is made explicit in \cite{RZnote}, for a particular situation considered there, that the inclusions in Lemma~\ref{lemarith} can be sharpened as follows:
$$\Phi_n^{-1}d_n^{s+1-i} \rho_i \in\Z \quad\mbox{for}\; i = 3, 5, \ldots, s,$$
and
$$\Phi_n^{-1}d_{n+1}^{s+1} \rho_{0,j} \in\Z \quad\mbox{for any}\; j\in \{1,\ldots,D\}, $$
where $\Phi_n=\Phi_n(D)$ is a certain product over primes in the range $2\le p\le n$, whose asymptotic behavior
$$
\phi=\phi(D)=\lim_{n\to\infty}\frac{\log\Phi_n}n
$$
can be controlled by means of the prime number theorem. It is possible to show that the quantity $\phi(D)/D$ increases to $\infty$ and at the same time $\phi(D)/(D\log^\varepsilon D)\to0$ as $D\to\infty$, for any choice of $\varepsilon>0$.
Later, we choose $D$ such that $ D\log D <s $, implying that the arithmetic gain coming from the factors $\Phi_n^{-1}$ is asymptotically negligible as $s\to\infty$.
\end{Rem}

\section{Asymptotic estimates of the linear forms}\label{sec4}

The following lemma is proved along the same lines as \cite[Lemma~2.1]{Sprang} (see also \cite[Lemma~4]{Zudilintrick} and the second proof of \cite[Lemme 3]{BR}). The  difference is that here we only assume $\frac{s}{  D\log D} $ to be sufficiently large, whereas in \cite{Sprang} parameter $D$ is fixed and $ s\to\infty$. 

\begin{Lem} \label{lemasy}
Assume that 
\begin{equation} \label{eqhyp}
\frac{s}{  D\log D} \quad\mbox{is sufficiently large.} 
\end{equation}
Then we have
\begin{equation} \label{eqasyun}
\lim_{n\to\infty} r_{n,j}^{1/n} = g(x_0) <3^{-(s+1)} \quad\mbox{and}\quad \lim_{n\to\infty} \frac{r_{n,j'}}{r_{n,j}} = 1
\quad\mbox{for any}\; j,j'\in\{1,\ldots,D\},
\end{equation}
where
$$g(x) = D^{ 3D } \, \frac{ (x+3)^{3D} (x+1)^{ s+1 }}{ ( x+2)^{2(s+1)}}$$
and $x_0$ is the unique positive root of the polynomial 
$$(X+3)^D(X+1)^{s+1} - X^D (X+2)^{s+1}.$$ 
\end{Lem}

\bigskip

\begin{Proof}
For $j\in \{1,\ldots,D\}$ and $k\geq 0$, let 
$$c_{k,j} = R_n\Big( n+k+\frac{j}{D}\Big) 
= D^{ 3D n} \, n!^{s+1-3D} \, \frac{ \prod_{\ell=0}^{3Dn} ( k +\frac{j+\ell }{D})}{ \prod_{\ell=0}^{ n} (n+k+\ell+\frac{j}{D})^{s+1}},$$
so that 
$$r_{n,j} = \sum_{m=1}^\infty R_n\Big(m+\frac{j}{D}\Big) = \sum_{k=0}^\infty c_{k,j}$$
is a sum  of positive terms.
We have
\begin{equation} \label{eqquok}
\frac{c_{k+1,j}}{c_{k,j}} = \bigg( \prod_{\ell = 1}^D \frac{k+3n+\frac{j+\ell}{D}}{k+\frac{j+\ell-1}{D}}\bigg) \, \bigg( \frac{k+n+\frac{j}{D}}{k+2n+1+ \frac{j}{D}}\bigg)^{s+1}
\end{equation}
implying that, for any $j$, the quotient $\frac{c_{k+1,j}}{c_{k,j}}$ tends to $ f(\kappa)$ as $n\to\infty$ assuming $k\sim \kappa n$ for $\kappa>0$ fixed, where 
$$f(x) = \Big( \frac{ x+ 3 }{x}\Big) ^D \Big( \frac{ x+ 1 }{x+2}\Big) ^{s+1}.$$
For the logarithmic derivative of this function we have
$$\frac{f'(x)}{f(x)} = 
\frac{D}{x+3} - \frac{D}{x }+\frac{s+1}{x+ 1 } - \frac{s+1}{x+ 2} =
\frac{ ax^2+bx+c}{x(x+1)(x+2)(x+3)}$$
with 
$a = s+1-3D > 0$ and $c = -6   D < 0$, hence the derivative $f'(x)$ vanishes exactly at one positive real number $x_1$. 
This means that the function $f(x)$ decreases on $(0, x_1]$ and increases on $[x_1,+\infty)$. Since $\lim_{x\to 0^+} f(x) = +\infty$ and $\lim_{x\to+\infty} f(x) = 1$, we deduce that there exists a unique positive real number $x_0$ such that $f(x_0 ) = 1$.

\bigskip

Let us now prove \eqref{eqasyun}. 
As in \cite[\S 3.4]{Bruijn} we wish to demonstrate that the asymptotic behaviour of $r_{n,j}$ is governed by the terms $c_{k,j}$ with $k$ close to $x_0n $ (see Eq.~\eqref{eqenca} below). To begin with, notice that 
\begin{align*}
c_{k,j} 
&= D^{-1} \, n!^{s+1-3D} \, \frac{ \prod_{\ell=0}^{3Dn} ( Dk+j+\ell )}{ \prod_{\ell=0}^{ n} (n+k+\ell+\frac{j}{D})^{s+1}}\\
&= D^{-1} \, n!^{s+1-3D} \, \frac{(3Dn+Dk+j)!}{(Dk+j-1)!} \, \frac{\Gamma(n+k+\frac{j}{D})^{s+1}}{\Gamma (2n+k +1+\frac{j}{D})^{s+1}}.
\end{align*}
Denoting by $k_0 (n) $ the integer part of $x_0 n$ and applying the Stirling formula to the factorial and gamma factors we obtain, as $n\to\infty$,
\begin{align}
c_{k_0(n),j}^{1/n} 
&\sim
\Big(\frac{n}{e}\Big) ^{s+1-3D} \,  \bigg(\frac{3Dn+Dk_0(n)+j }{e}\bigg) ^{ 3D+Dx_0} 
 \, \bigg(\frac{e}{Dk_0(n)+j-1 }\bigg) ^{ Dx_0} \nonumber \\
& \qquad \times 
\bigg(\frac{n+k_0(n)+\frac{j}{D} -1 }{e}\bigg) ^{ (s+1)( x_0+1)} 
 \, \bigg(\frac{e}{2n+k_0(n) +\frac{j}{D}}\bigg) ^{ (s+1)( x_0+2)} \nonumber \\
 &\sim \frac{((x_0+3)D)^{ (x_0+3)D}}{(x_0D)^{x_0D}}\, \frac{(x_0+1)^{(s+1)(x_0+1)} }{(x_0+2)^{(s+1)( x_0+2)} } \nonumber \\
& = g(x_0) f(x_0)^{x_0} = g(x_0). \label{eqasykz}
\end{align}

 \bigskip
We shall now give details that the asymptotic behavior of $r_{n,j}$ as $n\to\infty$ is determined by the terms $c_{k,j}$ with $k$ close to $x_0n $. Given $D$ and $s$, we take $\eps>0$ sufficiently small to accommodate the condition
$$b(\eps) = \max\Big( f(x_0+\eps), \frac1{f(x_0-\eps)}\Big)<1.$$
Then there exists $A(\eps) > x_1$, where $x_1$ is the unique positive root of $f'(x)=0$, such that $f(A(\eps))=b(\eps)$.
We have $f(x) \geq \frac1{b(\eps)}$ for any $x\in (0, x_0-\eps]$ and $f(x) \leq b(\eps) $ for any $x\in [ x_0+\eps, A(\eps)]$. For any $k$ such that $(x_0+2\eps)n \leq k\leq (A(\eps)-\eps)n$, Eq.~\eqref{eqquok} implies that $c_{k,j} \leq b(\eps) c_{k-1,j}$ provided $n$ is large (in terms of  $D$, $s$ and $\eps$), so that taking $k_1 = \lfloor (x_0+2\eps)n \rfloor$ and $k_2 = \lfloor (x_0+3\eps)n \rfloor$ we obtain
\begin{equation} \label{eqsommeun}
\sum_{k_2\leq k\leq (A(\eps)-\eps)n} c_{k,j} \leq c_{k_1,j}\sum_{k=k_2}^{+\infty} b(\eps)^{k-k_1}\leq c_{k_1,j}\frac{b(\eps)^{k_2-k_1}}{1-b(\eps)} \leq \eps \, c_{k_1,j}
\end{equation}
for all $n$ sufficiently large. In the same way, we get the estimate
\begin{equation} \label{eqsommede}
\sum_{1\leq k\leq \lfloor (x_0-3\eps)n \rfloor } c_{k,j} \leq \eps c_{ \lfloor (x_0-2\eps)n \rfloor , j}
\end{equation}
for all $n$ large (in terms of  $D$, $s$ and $\eps$).
At last, choosing $\eps$ small we can assume that $A(\eps)$ is sufficiently large (in terms of $D$ and $s$),
so that for $k \geq (A(\eps)-\eps)n$ we have
$$c_{k,j} \leq (2D)^{3Dn} \bigg( \frac{n!}{k^{n+1}}\bigg)^{ s+1 - 3D}$$
for $n$ large.
Using hypothesis \eqref{eqhyp} and the Stirling formula, the latter estimate implies
\begin{align}
\sum_{k = \lceil (A(\eps)-\eps)n\rceil} ^{+\infty} c_{k,j}
&\leq (3D)^{3Dn} \frac{n!^{ s+1 - 3D}}{((A(\eps)-\eps)n)^{(s+1-3D)(n+1)-1}} \nonumber\\
&\leq \bigg( \frac{3D}{e (A(\eps)-\eps)}\bigg)^{sn/2}
\leq \Big( \frac12 g(x_0)\Big)^n
\label{eqsommetr}
\end{align} 
provided $n$ is sufficiently large. Combining Eqs.~\eqref{eqasykz}, \eqref{eqsommeun}, \eqref{eqsommede} and \eqref{eqsommetr} we obtain
\begin{equation} \label{eqenca}
(1-3\eps) r_{n,j} \leq \sum_{ (x_0-3\eps)n \leq k \leq (x_0+3\eps)n } c_{k,j} \leq r_{n,j}.
\end{equation} 
Now for any $k$ in the range $ (x_0-3\eps)n \leq k \leq (x_0+3\eps)n $ it follows from the proof of Eq.~\eqref{eqasykz} that
$$g(x_0)-h(\eps) \leq c_{k,j}^{1/n} \leq g(x_0)+h(\eps)$$
for $n$ large (in terms of  $D$, $s$ and $\eps$), where $h$ is a positive function of $\eps$ such that $\lim_{\eps\to0^+} h(\eps) = 0$. This implies
$$(g(x_0)-2h(\eps))^n \leq 5\eps n (g(x_0)-h(\eps))^n \leq r_{n,j} \leq \frac{7\eps n}{1-3\eps} (g(x_0)+h(\eps))^n \leq (g(x_0)+2h(\eps))^n $$
for $n$ sufficiently large, and finishes the proof of $\lim_{n\to\infty} r_{n,j}^{1/n} = g(x_0)$ for any $j$.

\bigskip

To establish 
 $$\lim_{n\to\infty} \frac{r_{n,j'}}{r_{n,j}} = 1$$
for any $j,j'\in\{1,\ldots,D\}$, we can assume that $1\leq j \leq D-1$ and $j'=j+1$. For any $k $ we have
$$\frac{c_{k,j+1}}{c_{k,j}} = 
\frac{k+3D+\frac{j+1}{D}}{k+ \frac{j}{D}} \bigg( \frac{ \Gamma( n+k+\frac{j+1}{D})}{\Gamma( n+k +\frac{j }{D})}
\, \frac{ \Gamma( 2n+k+1+\frac{j }{D})}{\Gamma( 2n+k+1+\frac{j+1}{D})}\bigg)^{s+1}.$$
It follows from the Stirling formula that $\Gamma(x+\frac1{D}) \sim x^{1/D}\Gamma(x)$ as $x\to\infty$, so that for $k = \lfloor x_0 n\rfloor $ we have, as $n\to\infty$,
$$ \frac{c_{k,j+1}}{c_{k,j}} \sim
\frac{x_0+3}{x_0}\bigg( \frac{(x_0+1)^{1/D}}{(x_0+2)^{1/D}}\bigg)^{s+1} = f(x_0)^{1/D} = 1.$$
More generally, for $k$ in the range $ (x_0-3\eps)n \leq k \leq (x_0+3\eps)n $ and $n$ sufficiently large we have
$$1-\tilde h(\eps) \leq \frac{c_{k,j+1}}{c_{k,j}} \leq 1-\tilde h(\eps) $$
with $\lim_{\eps\to0^+} \tilde h(\eps) = 0$. Using Eq.~\eqref{eqenca} this concludes the proof of \eqref{eqasyun}, except for the upper bound on $g(x_0)$ that we shall verify now.

 \bigskip

To estimate $g(x_0)$ from above, we first show that $x_0<a$, where $a = 4\cdot 2^{-\frac{s+1}{D}}$.
Observe that $a < \frac12$, since $\frac{s}{D}\geq 3$. For any $x>0$ we have
$$ \frac{x+1}{x+2} \leq \frac{1}{2} \Big( 1 + \frac{x}{2}\Big)$$
implying
$$f(a) \leq \Big( \frac78\Big) ^D  \Big( 1 + \frac{a}{2}\Big)^{s+1}.$$
As $\frac sD$ is large and $\log(\frac78) < -\frac{1}{10}$, we deduce that
$$f(a)^{1/(s+1)} \leq \Big( \frac78\Big) ^{\frac{D}{s+1}} 
  \Big( 1 + \frac{a}{2}\Big) < \bigg( 1 -  \frac{1}{10} \frac{D}{s+1}\bigg)   \big( 1 +   2^{1-\frac{s+1}{D}}\big) < 1,$$
so that we indeed have $x_0<a<\frac12$.
Now this upper bound for $x_0$ implies
$$\log g(x_0) \leq 3D\log D + 3D  \log  \Big(\frac72\Big)  + (s+1) \big(  \log(a+1) - 2  \log (a+2)\big).$$
By taking $\frac{s}{D\log D}$ sufficiently large, we may ensure that the first two terms are sufficiently small in comparison with $s$ and that $a$ is sufficiently close to 0, so that $\log g(x_0) < - (s+1) \log 3$.

  This completes our proof of Lemma~\ref{lemasy}.
\end{Proof}

 \bigskip

\begin{Rem} \label{remexample}
For $s=77$ and $D=4$ one computes $g(x_0) <  \exp(-78)$. Thus, the suitable linear combinations
$$
		\hat{r}_{n,1}=r_{n,4}, \quad
		\hat{r}_{n,2}=r_{n,2}+r_{n,4} \quad\mbox{and}\quad
		\hat{r}_{n,4}=r_{n,1}+r_{n,2}+r_{n,3}+r_{n,4}
$$
of the corresponding linear forms allow us to eliminate three of the odd zeta values on the list
\[
\{\zeta(3),\zeta(5),\ldots,\zeta(77)\}.
\]
In particular, we obtain that two out of $\{\zeta(5),\zeta(7),\ldots,\zeta(77)\}$ are irrational. This result is slightly weaker than the result of Rivoal and the third author \cite{RZnote}, but it drops out as a byproduct of the construction above. The arithmetic gain given by $\Phi_n(4)$ for $\Phi_n(D)$ defined in Remark~\ref{remarith} can be used to slightly reduce the bound of $77$ to $73$, still weaker than the one in~\cite{RZnote}.
\end{Rem}

\section{A non-vanishing determinant}\label{sec5}

The following lemma~is used to eliminate  irrational zeta values in \S \ref{sec6} below. 

\begin{Lem} \label{lemdet}
For $t\geq 1$, let $x_1<\ldots<x_t$ be positive real numbers and $\alpha_1<\ldots<\alpha_t$ non-negative integers. Then the generalized Vandermonde matrix $[x_j^{\alpha_i}]_{1\leq i,j \leq t}$ has  positive determinant.
\end{Lem}

We remark that, subject to the hypothesis that $x_1,\ldots,x_t$ are real and positive, Lemma \ref{lemdet} is a stronger version of \cite[Lemme 1]{LMN} and, therefore, has potential applications to the zero estimates for linear forms in two logarithms.

\bigskip

The above result is quite classical and known to many people. While writing this paper we have found various proofs of rather different nature, three given below. We leave it to the readers to choose their favorite proof.

\begin{Proof}[Combinatorial proof of Lemma~\ref{lemdet}]
As pointed out in \cite[\S 2.1]{Krattdet}, the generalized Vandermonde determinant in question is closely related to Schur polynomials.
Let $\Delta := \det [x_j^{\alpha_{i}}]_{1\leq i,j \leq t}$, and 
$$V = \det [x_j^{i-1}]_{1\leq i,j \leq t} = \prod_{1\leq i < j \leq t} (x_j-x_i) >0$$
be the Vandermonde determinant of $x_1,\ldots, x_t$. For any $i\in\{1,\ldots,t\}$, we take $\lambda_i = \alpha_{t+1-i}+i - t $, so that $\lambda_1 \geq \ldots \geq \lambda_t \geq 0$; then $\lambda = (\lambda_1,\ldots,\lambda_t)$ is a partition of the integer $\lambda_1+\ldots+\lambda_t$. The associated Schur polynomial
$$s_\lambda = s_\lambda (x_1,\ldots, x_t)= \frac{\Delta}{V}$$
possesses the following expression:
$$s_\lambda = \sum_\mu K_{\lambda,\mu} m_\mu$$
with the sum over partitions $\mu$
(see, for instance, \cite[Appendix A]{FH} or \cite{Macdonald}, and \cite{Proctor} for a direct proof).
Here $m_\mu$ denotes the monomial symmetric polynomial $\sum_\sigma x_1^{\mu_{\sigma(1)}}\ldots  x_t^{\mu_{\sigma(t)}}$, where  the sum is over the distinct permutations of $\mu$,
while the coefficients $K_{\lambda,\mu} $
are non-negative integers and $K_{\lambda,\lambda}=1$. From this we deduce that $ s_\lambda$ is a positive real number, thus $\Delta=s_\lambda V>0$.
\end{Proof}

\begin{Proof}[Linear algebra proof of Lemma~\ref{lemdet}]
Write $A_{J,K}$ for the minor of an $n\times m$-matrix $A$, where $n\le m$, determined by ordered index sets $J$ and $K$.
A classical result due to Fekete \cite{fekete} asserts that if all $(n-1)$-minors
\[
A_{(1,2,\ldots,n-1),K},\quad K=(k_1,\ldots,k_{n-1}) \quad\text{with}\; 1\leq k_1<\ldots<k_{n-1}\leq m
\]
are positive, and all minors of size $n$ with consecutive columns are positive, then all $n$-minors of $A$ are positive. 
Thus, Lemma~\ref{lemdet} follows by induction on $t$ from Fekete's result applied to the matrix $[x_j^k]_{1\leq j \leq t,\, 0\leq k < m}$, using the positivity of the Vandermonde determinant.
\end{Proof}

\begin{Proof}[Analytical proof of Lemma~\ref{lemdet} {\rm(see \cite[p.~76--77]{GantmacherKrein})}]
By induction on $t$ one proves the following claim: {\em A non-zero function
\[
	f(x)=\sum_{i=1}^t c_i x^{\alpha_i}, 
\]
with $c_i,\alpha_i\in\R$, has  at most  $t-1$ positive zeros.}  Indeed, if $f$ has $t$ positive zeros then Rolle's theorem provides $t-1$ positive zeros of the derivative $\frac{\d}{\d x} (x^{-\alpha_1} f(x))$. The non-vanishing of the determinant in Lemma~\ref{lemdet} is an immediate consequence of this claim. Since the determinant depends continuously on the parameters $\alpha_i$, we deduce the required positivity from the positivity of the Vandermonde determinant.
\end{Proof}

\section{Elimination of odd zeta values}\label{sec6}

Let $0 < \eps < \frac13$, and let $s$ be odd and sufficiently large with respect to $\eps$. We take $D$ to be the product of all primes less than or equal to $(1-2\eps) \log s$ (such a product has asymptotically the largest possible number of divisors with respect to its size, see \cite[Chapter XVIII, \S 1]{HW}). We have
$$\log D \, = \sum_{\substack{p \; \mbox{\scriptsize prime}\\ p\leq (1-2\eps) \log s}} \log p \, \leq ( 1-\eps) \log s$$
by the prime number theorem, that is, $D\leq s^{1-\eps}$. Then $D \log D \leq  s^{1-\eps}\log s $: the assumption of Lemma~\ref{lemasy} holds. 

\bigskip

Notice that $D$ has precisely $\delta = 2^{\pi( (1-2\eps) \log s)}$ divisors, with 
$$\log \delta = \pi( (1-2\eps) \log s)\, \log 2 \geq (1-3\eps) (\log 2) \,\frac{\log s}{\log \log s}.$$
Assume that the number of irrational odd zeta values between $\zeta(3)$ and $\zeta(s)$ is less than $\delta$. Let $3 = i_1 < i_2 < \ldots < i_{\delta-1}\leq s$ be odd integers such that if $\zeta(i)\not\in\Q$ and $i$ is odd, $3\leq i \leq s$, then $i=i_j$ for some $j$. We set $i_0=1$, and consider the set $\calD$ of all divisors of $D$, so that $\Card \calD = \delta$. Lemma~\ref{lemdet} implies that the matrix $[d^{i_j}]_{d\in\calD, 0\leq j \leq \delta-1}$ is invertible. Therefore, there exist integers $w_d \in \Z$, where $d\in\calD$, such that 
\begin{equation} \label{eqwun}
\sum_{d\in\calD}w_d \, d^{i_j} = 0 \quad\mbox{for any}\; j\in\{1,\ldots,\delta-1\}
\end{equation}
and
\begin{equation} \label{eqwde}
\sum_{d\in\calD}w_d \, d^{i_0} = \sum_{d\in\calD}w_d \, d \neq 0.
\end{equation}
With the help of Lemma~\ref{lemconstru} we construct the linear forms 
$$r_{n,j} = \rho_{0,j}+\sum_{\substack{3\leq i \leq s\\i \;\mbox{\scriptsize odd}}} \rho_i \, \zeta\Big(i, \frac{j}{D}\Big)$$
for $n \geq 1$ and $1\leq j \leq D$. The crucial point (as in \cite[\S 3]{Sprang}) is that for any $d\in\calD$ and any $i\geq 2$,
$$\sum_{j=1}^d \zeta\bigg(i, \frac{j \frac{D}{d}}{D}\bigg) = \sum_{j=1}^d \zeta\Big(i, \frac{j}{d}\Big) = \sum_{n=0}^\infty \sum_{j=1}^d \frac{d^i}{(dn+j)^i} = d^i \zeta(i)$$
implying that
$$\widehat r_{n,d} = \sum_{j=1}^d r_{n, j\frac{D}{d}}
= \sum_{j=1}^d \rho_{0,j\frac{D}{d}} +\sum_{\substack{3\leq i \leq s\\i \;\mbox{\scriptsize odd}}} \rho_i \, d^i \, \zeta(i),$$
are linear forms in the odd zeta values with asymptotic behavior
$$\widehat r_{n,d} = (d+o(1))r_{n,1} \quad\text{as}\; n\to\infty, \quad\mbox{where}\; \lim_{n\to\infty} r_{n,1}^{1/n} = g(x_0) <3^{-(s+1)},$$
by Lemma~\ref{lemasy}.

\bigskip

We shall use now the integers $w_d$ to eliminate the odd zeta values $\zeta(i_j)$ for $j=1,\dots,\delta-1$, including all irrational ones, as in \cite{Zudilintrick} and \cite{Sprang}. For that, consider
$$\widetilde r_n = \sum_{d\in\calD} w_d \, \widehat r_{n,d}.$$
Eqs.~\eqref{eqwun} imply that
$$\widetilde r_n = \sum_{d\in\calD} w_d \sum_{j=1}^d \rho_{0,j\frac{D}{d}} +\sum_{i\in I } \rho_i \bigg( \sum_{d\in\calD} w_d \, d^i \bigg) \zeta(i),$$
where $I = \{3,5,7,\ldots,s\}\setminus \{i_1,\ldots,i_{\delta-1}\}$; in particular, no irrational zeta value $\zeta(i)$, where $3\leq i\leq s$, appears in this linear combination. Using Eq.~\eqref{eqwde} we obtain
$$ \widetilde r_n = \bigg( \sum_{d\in\calD}w_d \, d +o(1) \bigg) r_{n,1}  \quad \mbox{ with }   \quad \sum_{d\in\calD}w_d \, d \neq 0, $$
so that
$$\lim_{n\to\infty} | \widetilde r_n | ^{1/n} = g(x_0) <3^{-(s+1)}.$$

Now all $\zeta(i)$, $i\in I$, are assumed to be rational. Denoting by $A$ their common denominator, we deduce from Lemma~\ref{lemarith} that 
$A d_{n+1}^{s+1} \widetilde r_n $ is an integer. From the prime number theorem we have $\lim_{n\to\infty} d_{n+1}^{1/n} = e $, hence the sequence of integers satisfies 
$$0 < \lim_{n\to\infty} | A d_{n+1}^{s+1} \widetilde r_n | ^{1/n} = e^{s+1} g(x_0) < \Big(\frac{e}{3}\Big)^{s+1} < 1.$$ 
This contradiction concludes the proof of Theorem~\ref{th1}. 

\subsection*{Acknowledgements}
\noindent
We thank Michel Waldschmidt for his advice, Ole Warnaar for his comments on an earlier draft of the paper, and
Javier Fres\'an for educating us about the state of the art in Grothendieck's period conjecture and its consequences.

\newcommand{\url}{\texttt} \providecommand{\bysame}{\leavevmode ---\ }
\providecommand{\og}{``}
\providecommand{\fg}{''}
\providecommand{\smfandname}{\&}
\providecommand{\smfedsname}{\'eds.}
\providecommand{\smfedname}{\'ed.}
\providecommand{\smfmastersthesisname}{M\'emoire}
\providecommand{\smfphdthesisname}{Th\`ese}

\bigskip
\noindent
St\'ephane Fischler,
Laboratoire de Math\'ematiques d'Orsay, Univ. Paris-Sud, CNRS, Universit\'e Paris-Saclay, 91405 Orsay, France

\bigskip
\noindent
Johannes Sprang,
Fakult\"at f\"ur Mathematik, Universit\"at Regensburg, 93053 Regensburg, Germany

\bigskip
\noindent
Wadim Zudilin,
Department of Mathematics, IMAPP, Radboud University, PO Box 9010, 6500~GL Nijmegen, Netherlands; \\
School of Mathematical and Physical Sciences, The University of Newcastle, Callaghan, NSW 2308, Australia

\end{document}